\documentclass[12pt]{amsart}
\usepackage[cp1251]{inputenc}
\usepackage[T2A]{fontenc}
\usepackage[english]{babel}
\usepackage{amsmath,amsfonts,amssymb}
\usepackage{geometry}
\usepackage{amsmath, amsthm, amscd, amsfonts, amssymb, graphicx, color}
\usepackage[bookmarksnumbered, colorlinks, plainpages]{hyperref}

\newtheorem{theorem}{Theorem}

\newtheorem{proposition}{Proposition}

\theoremstyle{remark}

\theoremstyle{definition}

\begin{document}

\title[Sobolev-like Hilbert spaces induced by elliptic operators]{Sobolev-Like Hilbert Spaces\\Induced by Elliptic Operators}

\author[T. Kasirenko]{Tetiana Kasirenko}

\address{Institute of Mathematics, National Academy of Sciences of Ukraine,
3 Tereshchenkivs'ka, Kyiv, 01601, Ukraine}

\email{kasirenko@imath.kiev.ua}


\author[V. Mikhailets]{Vladimir Mikhailets}

\address{Institute of Mathematics, National Academy of Sciences of Ukraine,
3 Tereshchenkivs'ka, Kyiv, 01601, Ukraine}

\email{mikhailets@imath.kiev.ua}

%    Information for third author

\author[A. Murach]{Aleksandr Murach}

\address{Institute of Mathematics, National Academy of Sciences of Ukraine,
3 Tereshchenkivs'ka, Kyiv, 01601, Ukraine}

\email{murach@imath.kiev.ua}

%    General info



\begin{abstract}
We investigate properties of function spaces induced by the inner product Sobolev spaces $H^{s}(\Omega)$ over a bounded Euclidean domain $\Omega$ and by an elliptic differential operator $A$ on $\overline{\Omega}$. The domain and the coefficients of $A$ are of the class $C^{\infty}$. These spaces consist of all distributions $u\in H^{s}(\Omega)$ such that $Au\in H^{\lambda}(\Omega)$ and are endowed with the corresponding graph norm, with $s,\lambda\in\mathbb{R}$. We prove an interpolation formula for these spaces and discuss their application to elliptic boundary-value problems.
\end{abstract}

\maketitle

\section{Introduction}\label{sec1}

This paper is devoted to the function spaces induced by the scale $\{H^{s}(\Omega):s\in\mathbb{R}\}$ of inner product Sobolev spaces over a bounded Euclidean domain $\Omega$ and by an elliptic partial differential operator $A$ on $\overline{\Omega}$. The boundary of $\Omega$ and the coefficients of $A$ are supposed to be infinitely smooth. We investigate the inner product spaces $H^{s}_{A,\lambda}(\Omega)$, where $s,\lambda\in\mathbb{R}$, that consist of all distributions $u\in H^{s}(\Omega)$ subject to the condition $Au\in H^{\lambda}(\Omega)$ and that are endowed with the corresponding graph norm. We call $H^{s}_{A,\lambda}(\Omega)$ the Sobolev-like space induced by $A$. Such spaces are used in the theory of elliptic differential operators and elliptic boundary-value problems.

Specifically, the space $H^{0}_{A,0}(\Omega)$ is the domain of the maximal operator that corresponds to the unbounded operator $A$ defined on $C^{\infty}(\overline{\Omega})$ and acting in $L_{2}(\Omega)$ \cite{Ladyzhenskaya51, Hermander55, Birman57}. Lions and Magenes \cite{LionsMagenes62V, LionsMagenes63VI} systematically use the space $H^{s}_{A,0}(\Omega)$ for $s<\mathrm{ord}\,A=:2q$ in the theory of solvability of regular elliptic boundary-value problems in distribution spaces. In the case of the Dirichlet problem, the space $H^{s}_{A,-q}(\Omega)$ for $s<q$ is also useful \cite{LionsMagenes61II}.
These spaces serve as domains of Fredholm bounded operators that correspond to elliptic problems. A~more general class of the spaces $H^{s}_{A,\lambda}(\Omega)$ is considered in \cite{Murach09MFAT2, MikhailetsMurach14}. Some of their versions are introduced for $L_{p}$-Sobolev spaces \cite{LionsMagenes62V, LionsMagenes63VI, Geymonat62, Magenes65}, weighted Sobolev spaces \cite{Magenes65, LionsMagenes63VII, LionsMagenes72}, and generalized Sobolev spaces \cite{MikhailetsMurach14, 09OperatorTheory191, MikhailetsMurach12BJMA2} and are also used in the theory of elliptic differential equations.

The purpose of this paper is to investigate some properties of the spaces $H^{s}_{A,\lambda}(\Omega)$ concerning their completeness, separability, and dense subsets. Besides, we will prove an interpolation formula for these spaces. We also give their application to elliptic boundary-value problems.

\section{Main Results}\label{sec2}

Let $1\leq n\in\mathbb{Z}$ and $s\in\mathbb{R}$. We let $H^{s}(\mathbb{R}^{n})$ denote the complex inner product Sobolev space of order $s$ over $\mathbb{R}^{n}$. By definition, this space consists of all tempered distributions $w$ in $\mathbb{R}^{n}$ that their Fourier transform $\widehat{w}$ is locally Lebesgue integrable over $\mathbb{R}^{n}$ and satisfies the condition
$$
\|w\|_{s,\mathbb{R}^{n}}^{2}:=\int\limits_{\mathbb{R}^{n}}
(1+|\xi|^{2})^{s}\,|\widehat{w}(\xi)|^{2}d\xi<\infty.
$$
The space $H^{s}(\mathbb{R}^{n})$ is Hilbert and separable with respect to the norm $\|\cdot\|_{s,\mathbb{R}^{n}}$. If $\Omega$ is an open subset of $\mathbb{R}^{n}$, then the Sobolev space $H^{s}(\Omega)$ of order $s$ over $\Omega$ is defined to consist of the restrictions $u:=w\!\upharpoonright\!\Omega$ of all distributions $w\in H^{s}(\mathbb{R}^{n})$ to $\Omega$. This space is Hilbert and separable with respect to the norm
$$
\|u\|_{s,\Omega}:=
\inf\bigl\{\|w\|_{s,\mathbb{R}^{n}}:w\in
H^{s}(\mathbb{R}^{n}),\;u=w\!\upharpoonright\!\Omega\bigr\}
$$
and is continuously embedded in the linear topological space $\mathcal{D}'(\Omega)$ of all distributions in $\Omega$.

Henceforth we suppose that $n\geq2$ and that $\Omega$ is a bounded domain
whose boundary $\Gamma:=\partial\Omega$ is an infinitely smooth closed manifold of dimension $n-1$, the $C^\infty$-structure on $\Gamma$ being induced by $\mathbb{R}^{n}$. In the closed domain $\overline{\Omega}:=\Omega\cup\Gamma$, we consider an arbitrary elliptic linear partial differential expression
$$
A:=\sum_{|\mu|\leq 2q}a_{\mu}(x)\,
\frac{\partial^{|\mu|}}{\partial^{\mu_{1}}_{x_{1}}\cdots
\partial^{\mu_{n}}_{x_{n}}}
$$
of an even order $2q\geq2$ with complex-valued coefficients $a_{\mu}\in C^{\infty}(\overline{\Omega})$. Here,
$\mu:=(\mu_{1},\ldots,\mu_{n})$ is a multi-index with nonnegative integer-valued components, $|\mu|:=\mu_{1}+\cdots+\mu_{n}$, and $x=(x_1,\ldots,x_n)$ is an arbitrary point in $\mathbb{R}^{n}$. The ellipticity condition for $A$ means that
$$
\sum_{|\mu|=2q}a_{\mu}(x)\,
\xi_{1}^{\mu_{1}}\cdots\xi_{n}^{\mu_{n}}\neq0
$$
for every point $x\in\overline{\Omega}$ and each vector $\xi=(\xi_{1},\ldots,\xi_{n})\in\mathbb{R}^{n}\setminus\{0\}$. In the $n=2$ case, we suppose in addition that $A$ is properly elliptic on $\overline{\Omega}$ (see, e.g., \cite[Section~5.2.1]{Triebel95}).

We associate some spaces with the elliptic formal differential operator~$A$. Choosing $s,\lambda\in\mathbb{R}$ arbitrarily, we let $H^{s}_{A,\lambda}(\Omega)$ denote the linear space
$$
\bigl\{u\in H^{s}(\Omega):Au\in H^{\lambda}(\Omega)\bigr\}
$$
endowed with the graph norm
\begin{equation*}
\|u\|_{s,A,\lambda}:=
\bigl(\|u\|_{s,\Omega}^{2}+\|Au\|_{\lambda,\Omega}^{2}\bigr)^{1/2}.
\end{equation*}
Here and below, $Au$ is understood in the sense of the theory of distributions in $\Omega$.

\begin{theorem}\label{th1}
Let $s,\lambda\in\mathbb{R}$. The following assertions are true:
\begin{itemize}
\item[(i)] The space $H^{s}_{A,\lambda}(\Omega)$ is Hilbert and separable.
\item[(ii)] The equality of linear spaces $H^{s}_{A,\lambda}(\Omega)$ and $H^{s}(\Omega)$ holds if and only if $\lambda\leq s-2q$. If $\lambda\leq s-2q$, the norms in these spaces are equivalent.
\item[(iii)] The set $C^{\infty}(\overline{\Omega})$ is dense in $H^{s}_{A,\lambda}(\Omega)$.
\item[(iv)] If $s\leq1/2$ and $\lambda\leq1/2-2q$, then the set $C^{\infty}_{0}(\Omega)$ is dense in $H^{s}_{A,\lambda}(\Omega)$.
\end{itemize}
\end{theorem}

Here and below, $C^{\infty}_{0}(\Omega)$ denotes the set of all functions from $C^{\infty}(\overline{\Omega})$ that vanish near the boundary of $\Omega$.

Note that, if $\lambda>s-2q$, the space $H^{s}_{A,\lambda}(\Omega)$ depends on the coefficients of $A$, even when all of them are constant. For example, if $A_{1}$ and $A_{2}$ are elliptic constant-coefficient linear differential expressions, the embedding $H^{0}_{A_{1},0}(\Omega)\subseteq H^{0}_{A_{2},0}(\Omega)$ implies that $A_{2}=\alpha A_{1}+\beta$ for certain $\alpha,\beta\in\mathbb{C}$. This follows from H\"ormander's result \cite[Theorem~3.1]{Hermander55}.

It is well known that the Sobolev scale $\{H^{s}(\Omega):s\in\mathbb{R}\}$ has the following interpolation property: if $s_{0},s_{1}\in\mathbb{R}$ and $0<\theta<1$, then
\begin{equation}\label{interp}
\bigl[H^{s_{0}}(\Omega),H^{s_{1}}(\Omega)]_{\theta}=
H^{(1-\theta)s_{0}+\theta s_{1}}(\Omega)
\end{equation}
up to equivalence of norms (see, e.g., \cite[Section~4.3.1, Theorem~1]{Triebel95}). Here and below, $[H_{0},H_{1}]_{\theta}$ denotes the result of the complex interpolation with the parameter $\theta$ of a compatible pair of Hilbert (or, more generally, Banach) spaces $H_{0}$ and $H_{1}$ (see, e.g., \cite[Section~1.9]{Triebel95}). The pair $[H_{0},H_{1}]$ of these spaces is called compatible if they are continuously embedded in a certain linear Hausdorff space.

Let us formulate a version of the interpolation property \eqref{interp} for the spaces $H^{s_{j}}_{A,\lambda_{j}}(\Omega)$, where $j\in\{0,1\}$. In view of assertion~(ii) of Theorem~\ref{th1}, it is worthwhile to restrict ourselves to the case when $\lambda_{j}\geq s_{j}-2q$.

\begin{theorem}\label{th2}
Suppose that $s_{0}$, $s_{1}$, $\lambda_{0}$, $\lambda_{1}$, and $\theta$ are real numbers satisfying the inequalities $\lambda_{0}\geq s_{0}-2q$,  $\lambda_{1}\geq s_{1}-2q$, and $0<\theta<1$. Put $s:=(1-\theta)s_{0}+\theta s_{1}$ and $\lambda:=(1-\theta)\lambda_{0}+\theta\lambda_{1}$. Then
\begin{equation}\label{interp-th2}
\bigl[H^{s_{0}}_{A,\lambda_{0}}(\Omega),
H^{s_{1}}_{A,\lambda_{1}}(\Omega)]_{\theta}=
H^{s}_{A,\lambda}(\Omega)
\end{equation}
up to equivalence of norms.
\end{theorem}

\section{Proof of the main results}\label{sec3}

We will give a joint proof of Theorems \ref{th1} and \ref{th2}.

Let us first prove that the space $H^{s}_{A,\lambda}(\Omega)$ is Hilbert for arbitrary $s,\lambda\in\mathbb{R}$. Evidently, the norm in it is induced by a graph inner product. Besides, the space $H^{s}_{A,\lambda}(\Omega)$ is complete. Indeed, if $(u_{k})$ is a Cauchy sequence in this space, then there exist the limits $u:=\lim u_{k}$ in $H^{s}(\Omega)$ and $f:=\lim Au_{k}$ in $H^{\lambda}(\Omega)$ because the last two spaces are complete. The differential operator $A$ is continuous on $\mathcal{D}'(\Omega)$; therefore,
$Au=\lim Au_{k}=f$ in $\mathcal{D}'(\Omega)$. Here, $u\in H^{s}(\Omega)$ and $f\in H^{\lambda}(\Omega)$. Hence, $u\in H^{s}_{A,\lambda}(\Omega)$, and $\lim u_{k}=u$ in the space $H^{s}_{A,\lambda}(\Omega)$. Thus, this space is complete.

Let us now prove assertion~(ii) of Theorem~\ref{th1}. Since $\mathrm{ord}\,A=2q$, the differential operator $A$ acts continuously from $H^{s}(\Omega)$ to $H^{s-2q}(\Omega)$ for every $s\in\mathbb{R}$. Hence, if $\lambda\leq s-2q$, then $H^{s}(\Omega)=H^{s}_{A,\lambda}(\Omega)$ up to equivalence of norms due to the continuous embedding $H^{s-2q}(\Omega)\subseteq H^{\lambda}(\Omega)$.

Conversely, assume that $H^{s}(\Omega)=H^{s}_{A,\lambda}(\Omega)$ for certain $s,\lambda\in\mathbb{R}$. Choose an arbitrary distribution $w\in H^{s}(\mathbb{R}^{n})$ such that $\mathrm{supp}\,w\subset\Omega$, and put
$u:=w\!\upharpoonright\!\Omega$. Since $u\in H^{s}(\Omega)$, we have the inclusion $Au\in H^{\lambda}(\Omega)$ by our assumption. Hence, $\chi u\in H^{\lambda+2q}(\Omega)$ for every function $\chi\in C^{\infty}_{0}(\Omega)$ by virtue of the ellipticity of $A$ (see, e.g., \cite[Theorem 7.4.1]{Hermander63}). We can choose the function $\chi$ so that $\chi=1$ on $\mathrm{supp}\,w$. Therefore $u\in H^{\lambda+2q}(\Omega)$, which implies the inclusion $w\in H^{\lambda+2q}(\mathbb{R}^{n})$. Thus,
$$
\bigl\{w\in H^{s}(\mathbb{R}^{n}):\mathrm{supp}\,w\subset\Omega\bigr\}
 \subseteq H^{\lambda+2q}(\mathbb{R}^{n}).
$$
According to \cite[Theorem 2.2.2]{Hermander63}, this embedding implies that
$$
(1+|\xi|)^{\lambda+2q}\leq c\,(1+|\xi|)^{s}
$$
for every $\xi\in\mathbb{R}^{n}$ with some number $c>0$ not depending on $\xi$. Hence, $\lambda+2q\leq s$. We have proved assertion~(ii) of Theorem~\ref{th1}.

To prove Theorem~\ref{th2} and the remaining part of Theorem~\ref{th1}, we use a theorem on interpolation of some Hilbert spaces induced by linear continuous operators. Before we formulate this result, let us make some notation. If $H$, $\Phi$ and $\Psi$ are Hilbert spaces subject to the continuous embedding $\Phi\subseteq\Psi$ and if $T:H\rightarrow\Psi$ is a continuous linear operator, we put
$$
(H)_{T,\Phi}:=\{u\in H:\,Tu\in\Phi\}
$$
and endow the linear space $(H)_{T,\Phi}$ with the graph norm
$$
\|u\|_{(H)_{T,\Phi}}:=\bigl(\|u\|_{H}^{2}+\|Tu\|_{\Phi}^{2}\bigr)^{1/2}.
$$
This norm does not depend on $\Psi$, and the space $(H)_{T,\Phi}$ is Hilbert. The latter is proved in a quite similar way as the proof of the completeness of $H^{s}_{A,\lambda}(\Omega)$.

\begin{proposition}\label{prop-interp}
Assume that six separable Hilbert spaces $X_{0}$, $Y_{0}$, $Z_{0}$, $X_{1}$, $Y_{1}$, and $Z_{1}$ and three linear mappings $T$, $R$, and $S$ are given and satisfy the following conditions:
\begin{itemize}
\item[(i)] The pairs $[X_{0},X_{1}]$ and  $[Y_{0},Y_{1}]$ are compatible.
\item[(ii)] The spaces $Z_{0}$ and $Z_{1}$ are linear manifolds in a certain locally convex topological space $E$.
\item[(iii)] The continuous embeddings $Y_{0}\subseteq Z_{0}$ and $Y_{1}\subseteq Z_{1}$ hold.
\item[(iv)] The mapping $T$ is given on $X_{0}+X_{1}$ and defines the bounded operators $T:\nobreak X_{0}\rightarrow Z_{0}$ and $T:X_{1}\rightarrow Z_{1}$.
\item[(v)] The mapping $R$ is given on $E$ and defines the bounded operators $R:Z_{0}\rightarrow X_{0}$ and $R:Z_{1}\rightarrow X_{1}$.
\item[(vi)] The mapping $S$ is given on $E$ and defines the bounded operators $S:Z_{0}\rightarrow Y_{0}$ and
    $S:Z_{1}\rightarrow Y_{1}$.
\item[(vii)] The equality $TRu=u+Su$ holds for every $u\in E$.
\end{itemize}
Then
\begin{equation}\label{prop-interp-form}
[(X_{0})_{T,Y_{0}},(X_{1})_{T,Y_{1}}]_{\theta}=
([X_{0},X_{1}]_{\theta})_{T,[Y_{0},Y_{1}]_{\theta}}
\end{equation}
up to equivalence of norms for every $\theta\in(0,1)$.
\end{proposition}

This result is proved by Lions and Magenes \cite[Chapter~1, Theorem~14.3]{LionsMagenes72} in a more general case of Banach spaces.

Let us now prove Theorem~\ref{th2}. Our reasoning is motivated by \cite[Chapter~2, Proof of Theorem~4.2]{LionsMagenes72} and \cite[Section~3.3.4]{MikhailetsMurach14}. We choose an integer $r\geq1$ arbitrarily and consider the linear differential expression $A^{r}A^{r+}+{1}$ of order $4qr$. Here, as usual, $A^{r+}$ denotes the differential expression which is formally adjoint to the $r$-th iteration $A^{r}$ of~$A$. Given an integer $\sigma\geq2qr$, we let $H^{\sigma}_{D}(\Omega)$ denote the linear manifold of all distributions $u\in H^{\sigma}(\Omega)$ such that $\partial^{j}_{\nu}u=0$ on $\Gamma$ for each $j\in\{0,\ldots,2qr-1\}$. Here, $\partial_{\nu}$ is the operator of the differentiation with respect to the inward normal to the boundary $\Gamma$ of~$\Omega$. This definition is reasonable in view of the theorem on traces of distributions from $H^{\sigma}(\Omega)$ (see, e.g., \cite[Section~4.7.1]{Triebel95}). According to this theorem, $H^{\sigma}_{D}(\Omega)$ is a (closed) subspace of $H^{\sigma}(\Omega)$.
We have the isomorphism
\begin{equation*}
A^{r}A^{r+}+1:\;H^{\sigma}_{D}(\Omega)\leftrightarrow
H^{\sigma-4qr}(\Omega)
\end{equation*}
for every integer $\sigma\geq2qr$ (see, e.g., \cite[Lemma~3.1]{MikhailetsMurach14}). Let $(A^{r}A^{r+}+1)^{-1}$ denote the inverse of this isomorphism. This inverse defines the continuous linear operator
\begin{equation}\label{proof-inverse}
(A^{r}A^{r+}+1)^{-1}:H^{l}(\Omega)\rightarrow H^{l+4qr}(\Omega)
\end{equation}
for every integer $l\geq-2qr$. It follows from the interpolation formula
\eqref{interp} that this operator is well defined and continuous for every real $l\geq-2qr$.

In Proposition~\ref{prop-interp}, we put $X_{j}:=H^{s_{j}}(\Omega)$, $Y_{j}:=H^{\lambda_{j}}(\Omega)$, $Z_{j}:=H^{s_{j}-2q}(\Omega)$,  $E:=H^{\min\{s_{0},s_{1}\}-2q}(\Omega)$, and $T:=A$. Evidently, conditions (i)--(iv) of this proposition are fulfilled. Besides, subjecting $r$ to the restrictions
\begin{equation*}
s_{j}-2q\geq-2qr\quad\mbox{and}\quad s_{j}-2q-\lambda_{j}\geq-4qr
\end{equation*}
for each $j\in\{0,1\}$, we put $R:=A^{r-1}A^{r+}(A^{r}A^{r+}+1)^{-1}$ and $S:=-(A^{r}A^{r+}+1)^{-1}$. Owing to \eqref{proof-inverse}, we have the continuous linear operators
\begin{equation*}
R:Z_{j}=H^{s_{j}-2q}(\Omega)\to H^{s_{j}}(\Omega)=X_{j}
\end{equation*}
and
\begin{equation*}
S:Z_{j}=H^{s_{j}-2q}(\Omega)\to H^{s_{j}-2q+4qr}(\Omega)\subseteq
H^{\lambda_{j}}(\Omega)=Y_{j}
\end{equation*}
for each $j\in\{0,1\}$ (the last embedding is continuous); i.e., conditions (v) and (vi) are fulfilled as well. The last condition (vii) is satisfied because
\begin{equation*}
TRu=(A^{r}A^{r+}+1-1)(A^{r}A^{r+}+1)^{-1}u=u+Su
\end{equation*}
for every $u\in E$. Thus, owing to Proposition~\ref{prop-interp} and the interpolation formula \eqref{interp}, we conclude that
\begin{align*}
\bigl[H^{s_{0}}_{A,\lambda_{0}}(\Omega),
H^{s_{1}}_{A,\lambda_{1}}(\Omega)]_{\theta}&=
[(X_{0})_{T,Y_{0}},(X_{1})_{T,Y_{1}}]_{\theta}\\
&=([X_{0},X_{1}]_{\theta})_{T,[Y_{0},Y_{1}]_{\theta}}=
H^{s}_{A,\lambda}(\Omega)
\end{align*}
up to equivalence of norms, which proves Theorem~\ref{th2}.

Let us prove the remaining part of Theorem~\ref{th1}. We make use of the following supplement to Proposition~\ref{prop-interp}:

\begin{proposition}\label{prop2}
Suppose that the hypothesis of Proposition~$\ref{prop-interp}$ is fulfilled. If the Hilbert spaces $X_{0}$, $Y_{0}$, $X_{1}$, and $Y_{1}$ are separable and satisfy the dense continuous embeddings $X_{1}\subseteq X_{0}$ and $Y_{1}\subseteq Y_{0}$, then the Hilbert spaces $(X_{0})_{T,Y_{0}}$ and $(X_{1})_{T,Y_{1}}$ are also separable and satisfy the dense continuous embedding $(X_{1})_{T,Y_{1}}\subseteq(X_{0})_{T,Y_{0}}$.
\end{proposition}

This result is proved in \cite[Theorem~4.1]{MikhailetsMurach06UMJ11} for a more general method of interpolation between Hilbert spaces than that used in the paper (see also \cite[Section~3.3.2]{MikhailetsMurach14}).

According to Proposition~\ref{prop2}, the spaces $H^{s_{0}}_{A,\lambda_{0}}(\Omega)$ and $H^{s_{1}}_{A,\lambda_{1}}(\Omega)$ from Theorem~\ref{th2} are  separable, and the second of them is continuously and densely embedded in the first provided that $s_{0}\leq s_{1}$ and $\lambda_{0}\leq \lambda_{1}$. This specifically implies assertion (i) of Theorem~\ref{th1} concerning the separability of $H^{s}_{A,\lambda}(\Omega)$ for all $s\in\mathbb{R}$
and $\lambda>s-2q$. If $\lambda\leq s-2q$, this separability follows from assertion~(ii).

It remains to prove assertions (iii) and (iv) of Theorem~\ref{th1}. Let $s,\lambda\in\mathbb{R}$; if $\lambda>s-2q$, we have the dense continuous embedding
\begin{equation}\label{proof-embedding}
H^{s_{1}}(\Omega)=H^{s_{1}}_{A,s_{1}-2q}(\Omega)\subseteq
H^{s}_{A,\lambda}(\Omega)
\end{equation}
for every number $s_{1}\in\mathbb{R}$ such that $s_{1}\geq s$ and $s_{1}-2q\geq\lambda$. This has been noted in the previous paragraph. Since the set $C^{\infty}(\overline{\Omega})$ is dense in $H^{s_{1}}(\Omega)$, this embedding implies assertion~(iii) in the case of $\lambda>s-2q$. In the opposite case, this assertion is a consequence of assertion~(ii).

Finally, assume that $s\leq1/2$ and $\lambda\leq1/2-2q$. If $\lambda>s-2q$, assertion~(iv) follows from \eqref{proof-embedding} and the density of $C^{\infty}_{0}(\Omega)$ in $H^{s_{1}}(\Omega)$ for $s_{1}\leq1/2$ (see, e.g., \cite[Section~4.3.2, Theorem~1~(a)]{Triebel95}). If $\lambda\leq s-2q$, this assertion follows from this density in view of assertion~(ii). Our proof is complete.

\section{Application}\label{sec4}

Let us discuss an application of the space $H^{s}_{A,\lambda}(\Omega)$ to elliptic boundary-value problems. We consider a boundary-value problem that consists of the elliptic equation
\begin{equation}\label{1f1}
Au=f\quad\mbox{in}\quad\Omega
\end{equation}
and the boundary conditions
\begin{equation}\label{1f2}
B_{j}u=g_{j}\quad\mbox{on}\quad\Gamma,
\quad j=1,...,q.
\end{equation}
Here, each
$$
B_{j}:=\sum_{|\mu|\leq m_{j}}b_{j,\mu}(x)
\,\frac{\partial^{|\mu|}}{\partial^{\mu_{1}}_{x_{1}}\cdots
\partial^{\mu_{n}}_{x_{n}}}
$$
is a linear boundary partial differential expression on $\Gamma$ of order $m_{j}\leq2q-1$ with complex-valued coefficients $b_{j,\mu}\in C^{\infty}(\Gamma)$. We suppose that the boundary-value problem \eqref{1f1}, \eqref{1f2} is elliptic in $\Omega$ (see the definition in, e.g., \cite[Section~1.2]{Agranovich97}). Put $m:=\max\{m_{1},\ldots,m_{q}\}$; the case of $m\geq2q$ is possible.

It is known \cite[Theorem~7]{Peetre61} that the mapping
\begin{equation}\label{1f3}
\begin{gathered}
u\mapsto(Au,Bu)=(Au,B_{1}u,...,B_{q}u),
\end{gathered}
\end{equation}
where $u\in C^{\infty}(\overline{\Omega})$, extends uniquely (by continuity) to a Fredholm continuous linear operator
\begin{equation*}
(A,B):H^{s}(\Omega)\to H^{s-2q}(\Omega)\oplus
\bigoplus_{j=1}^{q}H^{s-m_{j}-1/2}(\Gamma)
\end{equation*}
for every real number $s>m+1/2$. Moreover, the kernel of this operator lies in $C^{\infty}(\overline{\Omega})$ and together with the index of the operator does not depend on $s$. Here and below, $H^{\sigma}(\Gamma)$ denotes the inner product Sobolev space over $\Gamma$ of order $\sigma\in\mathbb{R}$. Let $N$ and $\varkappa$ respectively denote these kernel and index.

This fundamental property of problem \eqref{1f1}, \eqref{1f2} does not remain for $s\leq m+1/2$. Specifically, if $s\leq m_{j}+1/2$, then the mapping $u\mapsto B_{j}u$, where $u\in C^{\infty}(\overline{\Omega})$, cannot be extended to a continuous linear operator from $H^{s}(\Omega)$ to the linear topological space $\mathcal{D}'(\Gamma)$ of all distributions on $\Gamma$. But it is possible to prove a version of this property for every $s\leq m+1/2$ provided that we use some spaces $H^{s}_{A,\lambda}(\Omega)$ instead of $H^{s}(\Omega)$.

\begin{theorem}\label{th3}
Let numbers $s,\lambda\in\mathbb{R}$ satisfy the conditions $s\leq m+1/2$, $\lambda>-1/2$, and $\lambda>m+1/2-2q$. Then the  mapping \eqref{1f3}, where $u\in C^{\infty}(\overline{\Omega})$, extends uniquely (by continuity) to a Fredholm continuous linear operator
\begin{equation}\label{3f14}
(A,B):H^{s}_{A,\lambda}(\Omega)\to
H^{\lambda}(\Omega)\oplus
\bigoplus_{j=1}^{q}H^{s-m_{j}-1/2}(\Gamma).
\end{equation}
The kernel of this operator equals $N$, and the index equals $\varkappa$.
\end{theorem}

If $s\in\mathbb{Z}$ and $s\leq0$, Theorem~\ref{th3} is proved in a somewhat similar way as that used in \cite{Murach09MFAT2} and \cite[Sections 4.4.2 and 4.4.3]{MikhailetsMurach14} for regular elliptic boundary-value problems. (The corresponding reasoning is given in \cite[Section~5]{AnopKasirenkoMurach18UMJ3} for $m\geq2q$). In the general situation we prove Theorem~\ref{th3} with the help of the interpolation formula \eqref{interp-th2}. Choose $l\in\mathbb{Z}$ such that $l<s$ and $l\leq0$. We consider the Fredholm continuous operators
\begin{equation}\label{operator-1}
(A,B):H^{l}_{A,\lambda}(\Omega)\to H^{\lambda}(\Omega)\oplus
\bigoplus_{j=1}^{q}H^{l-m_{j}-1/2}(\Gamma)
\end{equation}
and
\begin{equation}\label{operator-2}
(A,B):H^{\lambda+2q}(\Omega)\to H^{\lambda}(\Omega)\oplus
\bigoplus_{j=1}^{q}H^{\lambda+2q-m_{j}-1/2}(\Gamma)
\end{equation}
(note that $\lambda+2q>1/2$ by the hypothesis of the Theorem). They have the common kernel $N$ and the common index $\varkappa$. Since $l<s<\lambda+2q$, there exists a number $\theta\in(0,1)$ such that
$s=(1-\theta)l+\theta(\lambda+2q)$. Applying the interpolation with the parameter $\theta$ to these operators, we conclude by
\cite[Theorem~1.5]{MikhailetsMurach14} that a restriction of the first operator is a Fredholm continuous operator between the spaces
\begin{equation*}
\bigl[H^{l}_{A,\lambda}(\Omega),H^{\lambda+2q}(\Omega)]_{\theta}
\end{equation*}
and
\begin{equation*}
\biggl[H^{\lambda}(\Omega)\oplus
\bigoplus_{j=1}^{q}H^{l-m_{j}-1/2}(\Gamma),
H^{\lambda}(\Omega)\oplus
\bigoplus_{j=1}^{q}H^{\lambda+2q-m_{j}-1/2}(\Gamma)\biggr]_{\theta}.
\end{equation*}
Moreover, this operator has the same kernel and index as operators \eqref{operator-1} and \eqref{operator-2}. Owing to Theorem~\ref{th2}, the first space coincides with $H^{s}_{A,\lambda}(\Omega)$ up to equivalence of norms (remark that $H^{\lambda+2q}(\Omega)=H^{\lambda+2q}_{A,\lambda}(\Omega)$). Besides, the second space coincides up to equivalence of norms with the target space of operator \eqref{3f14} due to an analog of the interpolation formula \eqref{interp} for Sobolev spaces over $\Gamma$ (see, e.g, \cite[Capter~1, Theorem~7.7]{LionsMagenes72}). Thus, the latter Fredholm operator acts between the spaces indicated in \eqref{3f14}. It is an extension by continuity of the mapping \eqref{1f3}, where $u\in C^{\infty}(\overline{\Omega})$, in view of assertion (iii) of Theorem~\ref{th1}. The proof of Theorem~\ref{th3} is complete.

Remark that, in the important case of regular elliptic boundary-value problems, Theorem~\ref{th3} is proved by Lions and Magenes \cite{LionsMagenes62V, LionsMagenes63VI} in the framework of $L_{p}$-Sobolev spaces provided that $\lambda=0$ and $s\notin\{1/p-k:1\leq k\in\mathbb{Z}\}$. In the case where $s\geq0$ and $m\leq2q-1$, this theorem is contained in the result formulated in \cite[p.~86]{Agranovich97}.


\begin{thebibliography}{99}

\bibitem{Ladyzhenskaya51}
O. A. Ladyzenskaya, \textit{On the closure of the elliptic operator} (Russian), Doklady Akad. Nauk SSSR (N.S.) \textbf{79}, (1951). 723–-725.

\bibitem{Hermander55}
L. H\"ormander, \textit{On the theory of general partial differential
equations}, Acta Math. \textbf{94} (1955), no.~1, 161--248.

\bibitem{Birman57}
M. \v{S}. Birman, \textit{Characterization of elliptic differential operators with the maximal domain of definition} (Russian), Vestnik Leningrad. Univ. Ser. Math. Mech. Astr. \textbf{19}  (1957), no.~4, 177--183.

\bibitem{LionsMagenes62V}
J.-L. Lions, E. Magenes, \textit{Probl\'emes aux limites  non homog\'enes, V}, Ann. Sci. Norm. Sup. Pisa \textbf{16} (1962), 1--44.

\bibitem{LionsMagenes63VI}
J.-L. Lions, E. Magenes, \textit{Probl\'emes aux limites non homog\'enes, VI}, J. d'Analyse Math. \textbf{11} (1963), 165--188.

\bibitem{LionsMagenes61II}
J.-L. Lions, E. Magenes, \textit{Probl\'emes aux limites  non homog\'enes, II}, Ann. Inst. Fourier (Grenoble) \textbf{11} (1961), 137--178.

\bibitem{Murach09MFAT2}
Murach A. A., \textit{Extension of some Lions--Magenes theorems},
Methods Funct. Anal. Topology, \textbf{15} (2009), no.~2, 152--167.

\bibitem{MikhailetsMurach14}
V. A. Mikhailets, A. A. Murach, \textit{H\"ormander spaces, interpolation, and elliptic problems}, de Gruyter Stud. Math., vol.~60, De Gruyter, Berlin, 2014.

\bibitem{Geymonat62}
G. Geymonat, \textit{Sul problema di Dirichlet per le equazoni lineari ellittiche}, Ann. Sci. Norm. Sup. Pisa \textbf{16} (1962), 225--284.

\bibitem{Magenes65}
E. Magenes, \textit{Spazi di interpolazione ed equazioni a derivate parziali}, Atti VII Congr. Un. Mat. Italiana (Genoa, 1963), Edizioni Cremonese, Rome, 1965, 134--197.

\bibitem{LionsMagenes63VII}
J.-L. Lions, E. Magenes, \textit{Probl\'emes aux limites  non homog\'enes, VII}, Ann. Mat. Pura Appl. (4) \textbf{63} (1963), 201--224.

\bibitem{LionsMagenes72}
J.-L. Lions and E. Magenes, \textit{Non-Homogeneous Boundary-Value Problems and Applications, vol.~I}, Grundlehren Math. Wiss., vol.~181,
Springer-Verlag, New York\,--\,Heidelberg, 1972.

\bibitem{09OperatorTheory191}
V. A. Mikhailets and A. A. Murach, \textit{Elliptic problems and H\"ormander spaces}, Modern analysis and applications. The Mark Krein Centenary Conference. Vol.~2: Differential operators and mechanics, 447--470, Oper. Theory Adv. Appl. \textbf{191}, Birkh\"aser, Basel, 2009.

\bibitem{MikhailetsMurach12BJMA2}
V. A. Mikhailets, A. A. Murach, \textit{The refined Sobolev scale,
inter\-po\-la\-tion, and elliptic problems}, Banach J. Math. Anal. \textbf{6} (2012), no.~2, 211--281.

\bibitem{Triebel95}
H.~Triebel, \textit{Interpolation Theory, Function Spaces, Differential Operators} [2-nd edn], Johann Ambrosius Barth, Heidelberg, 1995.

\bibitem{Hermander63}
L. H\"ormander, \textit{Linear Partial Differential Operators}, Grundlehren Math. Wiss., vol.~116, Springer, Berlin, 1963.

\bibitem{MikhailetsMurach06UMJ11}
V. A. Mikhailets, A. A. Murach, \textit{A regular elliptic boundary-value problem for a homogeneous equation in a two-sided refined scale of spaces}, Ukrainian Math.~J. \textbf{58} (2006), no.~11, 1748--1767.

\bibitem{Peetre61}
J. Peetre, \textit{Another approach to elliptiv boundary-value problems}, Comm. Pure Appl. Math. \textbf{14} (1961), no.~4, 711--731.

\bibitem{Agranovich97}
M. S. Agranovich, \textit{Elliptic boundary problems}, Partial differential equations, IX, 1--144, Encyclopaedia Math. Sci., \textbf{79}, Springer, Berlin, 1997.

\bibitem{AnopKasirenkoMurach18UMJ3}
A. V. Anop, T. M. Kasirenko, A. A. Murach, \textit{Nonregular elliptic boundary-value problems and H\"ormander spaces}, Ukra\"{\i}n. Mat. Zh. \textbf{70} (2018), no.~3, 299--317.

\end{thebibliography}
\end{document}